\documentclass[12pt,reqno]{article}

\usepackage[usenames]{color}
\usepackage{amssymb}
\usepackage{graphicx}
\usepackage{amscd}

\usepackage[colorlinks=true,
linkcolor=webgreen,
filecolor=webbrown,
citecolor=webgreen]{hyperref}

\definecolor{webgreen}{rgb}{0,.5,0}
\definecolor{webbrown}{rgb}{.6,0,0}

\usepackage{color}
\usepackage{fullpage}
\usepackage{float}

\usepackage{fullpage}
\usepackage{float}
\usepackage{graphics,amsmath,amssymb}
\usepackage{amsthm}
\usepackage{amsfonts}
\usepackage{latexsym}
\usepackage{epsf}
\usepackage{enumerate}

\newcommand{\seqnum}[1]{\href{http://oeis.org/#1}{\underline{#1}}}

\newcommand{\Z}{\mathbb Z}

\newcommand{\Q}{\mathbb Q}

\begin{document}

\theoremstyle{plain}
\newtheorem{theorem}{Theorem}
\newtheorem{corollary}[theorem]{Corollary}
\newtheorem{lemma}[theorem]{Lemma}
\newtheorem{proposition}[theorem]{Proposition}

\theoremstyle{definition}
\newtheorem{definition}[theorem]{Definition}
\newtheorem{example}[theorem]{Example}
\newtheorem{conjecture}[theorem]{Conjecture}

\theoremstyle{remark}
\newtheorem{remark}[theorem]{Remark}

\begin{center}
\vskip 1cm{\LARGE\bf 
On the Largest Integer that is not a Sum of Distinct Positive \(n\)th Powers
}
\vskip 1cm
\large
Doyon Kim\\
Department of Mathematics and Statistics\\
Auburn University\\
Auburn, AL 36849 \\
USA \\
\href{mailto:dzk0028@auburn.edu}{\tt dzk0028@auburn.edu} \\
\end{center}

\vskip .2 in
\begin{abstract}
It is known that for an arbitrary positive integer \(n\) the sequence \(S(x^n)=(1^n, 2^n, \ldots)\) is complete, meaning that every sufficiently large integer is a sum of distinct \(n\)th powers of positive integers. We prove that every integer \(m\geq (b-1)2^{n-1}(r+\frac{2}{3}(b-1)(2^{2n}-1)+2(b-2))^n-2a+ab\), where \(a=n!2^{n^2}\), \(b=2^{n^3}a^{n-1}\), \(r=2^{n^2-n}a\), is a sum of distinct positive \(n\)th powers.
\end{abstract}

\section{Introduction}
Let \(S=(s_1,s_2,\ldots)\) be a sequence of integers. The sequence \(S\) is said to be \textit{complete} if every sufficiently large integer can be represented as a sum of distinct elements of \(S\). For a complete sequence \(S\), the largest integer that is not representable as a sum of distinct elements of \(S\) is called the \textit{threshold of completeness} of \(S\). We let \(\theta_S\) denote the threshold of completeness of \(S\). 
\par
The threshold of completeness is often very difficult to find even for a simple sequence. For an arbitrary positive integer \(n\), let \(S(x^n)\) denote the sequence of \(n\)th powers of positive integers, i.e., \(S(x^n)=(1^n, 2^n, \ldots)\). The completeness of the sequence was proved in 1948, by Sprague \cite{Sprague6}. In 1954, Roth and Szekeres \cite{Roth5} further generalized the result by proving that if \(f(x)\) is a polynomial that maps integers into integers, then \(S(f)=(f(1),f(2), \ldots)\) is complete if and only if \(f(x)\) has a positive leading coefficient and for any prime \(p\) there exists an integer \(m\) such that \(p\) does not divide \(f(m)\). In 1964, Graham \cite{Graham2} re-proved the theorem of Roth and Szekeres using alternative elementary techniques.
\par
However, little is known about the threshold of completeness of \(S(x^n)\). The value \(\theta_{S(x^n)}\) is known only for \(n\leq 6\). The values are as follows: \(\theta_{S(x)}=0\), \(\theta_{S(x^2)}=128\) \cite{Sprague7}, \(\theta_{S(x^3)}=12758\) \cite{Graham2}, \(\theta_{S(x^4)}=5134240\) \cite{Lin3}, \(\theta_{S(x^5)}=67898771\) \cite{Patterson4}, \(\theta_{S(x^6)}=11146309947\) \cite{Fuller1}. Sprague, Roth and Szekeres, and Graham proved that \(S(x^n)\) is complete, but they were not interested in the size of \(\theta_{S(x^n)}\). The values \(\theta_{S(x^n)}\) for \(3\leq n\leq 6\) were found by methods that require lengthy calculations assisted by computer, and they do not give any idea on the size of \(\theta_{S(x^n)}\) for general \(n\).
\par
In this paper, we establish an upper bound of \(\theta_{S(x^n)}\) as a function of \(n\). Using the elementary techniques Graham used in his proof, it is possible to obtain an explicit upper bound of the threshold of completeness of \(S(x^n)=(1^n,2^n,3^n,\ldots)\). Since the case \(n=1\) is trivial, we let \(n\) be a positive integer greater than \(1\). We prove the following theorem:
\begin{theorem}
Let \(a=n!2^{n^2}\), \(b=2^{n^3}a^{n-1}\) and \(r=2^{n^2-n}a\). Then
\[\theta_{S(x^n)}<(b-1)2^{n-1}(r+\frac{2}{3}(b-1)(2^{2n}-1)+2(b-2))^n-2a+ab.\]
\end{theorem}
The theorem yields the result 
\[
\theta_{S(x^n)}=O((n!)^{n^2-1}\cdot 2^{2n^4+n^3+n^2+(2-\frac{\ln 3}{\ln 2})n}).
\]
The upper bound of \(\theta_{S(x^n)}\) given by the formula is much greater than \(4^{n^4}\), while the actual values of \(\theta_{S(x^n)}\) for \(2\leq n \leq 6\) are less than \(4^{n^2}\). So the upper bound obtained in this paper is most likely far from being tight.
\section{Preliminary results}
Let \(S=(s_1, s_2, \ldots)\) be a sequence of integers. \par
\begin{definition} The set \(P(S)\) is a set of all sums of the form \(\sum_{k=1}^{\infty} \epsilon_k s_k\) where \(\epsilon_k\) is \(0\) or \(1\), all but a finite number of \(\epsilon_k\) are \(0\) and at least one of \(\epsilon_k\) is \(1\).
\end{definition}
\begin{definition}
The sequence \(S\) is \textit{complete} if \(P(S)\) contains every sufficiently large integer.
\end{definition} 
\begin{definition}
If \(S\) is complete, the \textit{threshold of completeness} \(\theta_S\) is the largest integer that is not in \(P(S)\).
\end{definition}
\begin{definition}
The set \(A(S)\) is a set of all sums of the form \(\sum_{k=1}^{\infty} \delta_k s_k\) where \(\delta_k\) is \(-1\), \(0\) or \(1\) and all but a finite number of \(\delta_k\) are \(0\).
\end{definition}
\begin{definition}
Let \(k\) be a positive integer. The sequence \(S\) is a \textit{\(\Sigma(k)\)-sequence} if \(s_1\leq k\), and \[s_n\leq k+\sum_{j=1}^{n-1}s_j, \quad n\geq 2.\]
\end{definition}
For example, if \(S=(2,4,8,16,\ldots)\) then \(S\) is a \(\Sigma(2)\)-sequence since \(2^n=2+\sum_{j=1}^{n-1}2^j\) for all \(n\geq 2\).
\begin{definition}
Let \(c\) and \(k\) be positive integers. The sequence \(S\) is \textit{\((c,k)\)-representable} if \(P(S)\) contains \(k\) consecutive integers \(c+j\), \(1\leq j \leq k\).
\end{definition}
For example, if \(S=(1,3,6,10,\ldots)\) is a sequence of triangle numbers then \(S\) is \((8,3)\)-representable since \(\{9,10,11\}\subset P(S)\).
\begin{definition}
For a positive integer \(m\), we define \(\Z_m(S)\) to be the sequence \((\alpha_1,\alpha_2,\ldots)\), where \(0\leq \alpha_i<m\) and \(s_i\equiv \alpha_i\pmod m\) for all \(i\).
\end{definition}
The two following lemmas, slightly modified from Lemma 1 and Lemma 2 in Graham's paper \cite{Graham2}, are used to obtain the upper bound.
\begin{lemma} \label{1}
For a positive integer \(k\), let \(S=(s_1,s_2,\ldots)\) be a strictly increasing \(\Sigma(k)\)-sequence of positive integers and let \(T=(t_1,t_2,\ldots)\) be \((c,k)\)-representable. Then \(U=(s_1,t_1,s_2,t_2,\ldots)\) is complete and \(\theta_U\leq c\).
\begin{proof}
It suffices to prove that every positive integer greater than \(c\) belongs to \(P(U)\). The proof proceeds by induction. Note that all the integers \(c+t\), \(1\leq t\leq k\) belong to \(P(T)\), and all the integers \(c+s_1+t\), \(1\leq t\leq k\) belong to \(P(U)\). If \(1\leq t\leq k\) then \[c+t\in P(T)\subset P(U),\] and if \(k+1\leq t\leq k+s_1\), then \(1\leq k-s_1+1\leq t-s_1\leq k\) and we have \[c+t=c+(t-s_1)+s_1\in P(U).\] Therefore all the integers \[c+t, \quad 1\leq t\leq k+s_1\] belong to \(P(U)\). Now, let \(n\geq 2\) and suppose that all the integers \[c+t,\quad 1\leq t\leq k+\sum_{j=1}^{n-1} s_j\] belong to \(P(U)\), and that for every such \(t\) there is a \(P(U)\) representation of \(c+t\) such that none of \(s_m\), \(m\geq n\) is in the sum. Since all the integers \(c+t+s_n\), \(1\leq t\leq k+\sum_{j=1}^{n-1} s_j\) belong to \(P(U)\) and \(c+1+s_n\leq c+1+k+\sum_{j=1}^{n-1} s_j\), all the integers \[c+t, \quad 1\leq t\leq k+\sum_{j=1}^{n} s_j\] belong to \(P(U)\). Since \(S\) is a strictly increasing sequence of positive integers, this completes the induction step and the proof of lemma.
\end{proof}
\end{lemma}
\begin{lemma} \label{2}
Let \(S=(s_1,s_2,\ldots)\) be a strictly increasing sequence of positive integers. If \(s_k\leq 2s_{k-1}\) for all \(k\geq 2\), then \(S\) is a \(\Sigma(s_1)\)-sequence.
\begin{proof} For \(k\geq 2\), we have
\begin{align*}
s_k &\leq 2s_{k-1}=s_{k-1}+s_{k-1} \\
&\leq s_{k-1}+2s_{k-2}=s_{k-1}+s_{k-2}+s_{k-2} \\
&\leq s_{k-1}+s_{k-2}+2s_{k-3}\leq\cdots \\
&\leq s_1+\sum_{j=1}^{k-1}s_j.
\end{align*}
Therefore, \(S\) is a \(\Sigma(s_1)\)-sequence.
\end{proof}
\end{lemma}
Lemma \ref{1} shows that if a sequence \(S\) can be partitioned into one \(\Sigma(k)\)-sequence and one \((c,k)\)-representable sequence then \(S\) is complete with \(\theta_S\leq c\). What we aim to do is to partition \(S(x^n)\) into two such sequences for some \(c\) and \(k\). \par
Let \(f(x)=x^n\) and let \(S(f)=(f(1),f(2), \ldots)\). Let \(a=n!2^{n^2}\) and \(r=2^{n^2-n}a\). Partition the elements of the sequence \(S(f)\) into four sets \(B_1\), \(B_2\), \(B_3\) and \(B_4\) defined by
\begin{align*}
    B_1&=\{f(\alpha a+\beta):\; 0\leq \alpha\leq 2^{n^2-n}-1, \; 1\leq \beta\leq 2^n\}, \\
B_2&=\{f(\alpha a+\beta):\; 0\leq \alpha\leq 2^{n^2-n}-1, \; 2^n+1\leq \beta\leq a,\; \alpha a+\beta<2^{n^2-n} a\},\\
B_3&=\{f(2^{n^2-n}a),f(2^{n^2-n}a+2),f(2^{n^2-n}a+4),\ldots\},\\
B_4&=\{f(2^{n^2-n}a+1),f(2^{n^2-n}a+3),f(2^{n^2-n}a+5),\ldots\},
\end{align*}
so that 
\[B_1\cup B_2 =\{f(1),f(2),\ldots,f(r-1)\}\]
and
\[B_3\cup B_4 =\{f(r),f(r+1),f(r+2),\ldots\}.\]
Let \(S\), \(T\), \(U\) and \(W\) be the strictly increasing sequences defined by
\begin{align*}
S&=(s_1,s_2,\ldots,s_{2^{n^2}}),\quad  s_j\in B_1, \\ 
T&=(t_1,t_2,\ldots), \quad t_j\in B_3,\\
U&=(u_1,u_2,\ldots), \quad u_j\in B_1\cup B_3,\\
W&=(w_1,w_2,\ldots), \quad w_j\in B_2\cup B_4.
\end{align*}
Then the sequences \(U\) and \(W\) partition the sequence \(S(f)\). First, using Lemma \ref{2}, we show that \(W\) is a \(\Sigma(a)\)-sequence.
\begin{lemma} \label{3}
For \(a=n!2^{n^2}\) and \(r=2^{n^2-n}a\), \[\frac{f(r+1)}{f(r-1)}<\frac{f(a+2^n+1)}{f(a)}<\frac{f(2^n+2)}{f(2^n+1)}\leq 2.\]
\begin{proof}
Re-write the inequalities as 
\[\left(1+\frac{2}{r-1}\right)^n<\left(1+\frac{2^n+1}{a}\right)^n<\left(1+\frac{1}{2^n+1}\right)^n\leq 2.\]
It is clear that 
\[\frac{r-1}{2}>\frac{a}{2^n+1}>2^n+1,\]
which proves the first two inequalities.
The proof of the third inequality 
\[\left(1+\frac{1}{2^n+1}\right)^n\leq 2 \iff 1\leq (2^{\frac{1}{n}}-1)(2^n+1)\]
is also straightforward.
\end{proof}
\end{lemma}
\begin{corollary} \label{cor1}
The sequence \(W\) is a \(\Sigma(a)\)-sequence.
\begin{proof}
Note that \(w_1=(2^n+1)^n\). For every \(k\geq 2\), \(\frac{w_k}{w_{k-1}}\) satisfies one of the following equalities:
\begin{align}
    \frac{w_k}{w_{k-1}}&=\frac{f(\alpha+1)}{f(\alpha)}, \quad \text{for} \quad \alpha\geq 2^n+1;\\
    \frac{w_k}{w_{k-1}}&=\frac{f(\beta a+2^n+1)}{f(\beta a)},\quad \text{for}\quad \beta\geq 1; \\
    \frac{w_k}{w_{k-1}}&=\frac{f(\gamma+2)}{f(\gamma)},\quad \text{for}\quad \gamma\geq r-1.
\end{align}
Also, for every \(\alpha\geq 2^n+1\), \(\beta\geq 1\) and \(\gamma\geq r-1\) we have
\begin{align*}
\frac{f(\alpha+1)}{f(\alpha)}&\leq\frac{f(2^n+2)}{f(2^n+1)},\\
\frac{f(\beta a+2^n+1)}{f(\beta a)}&\leq\frac{f(a+2^n+1)}{f(a)},\\
\frac{f(\gamma+2)}{f(\gamma)}&\leq\frac{f(r+1)}{f(r-1)}.
\end{align*}
Thus, by Lemma \ref{3}, \(\frac{w_k}{w_{k-1}}\leq 2\) for \(k\geq 2\), and therefore by Lemma \ref{2}, \(W\) is a \(\Sigma((2^n+1)^n)\)-sequence. To complete the proof, it remains to prove that \((2^n+1)^n<a\) for all \(n>1\). The inequality is true for \(n=2\) and \(n=3\), and for \(n>3\) we have
\[(2^n+1)^n<(2^n+2^n)^n=2^n2^{n^2}<n!2^{n^2}=a.\] 
Therefore, \(W\) is a \(\Sigma(a)\)-sequence.
\end{proof}
\end{corollary}
Now, we prove that \(U\) is \((d,a)\)-representable for some positive integer \(d\). By Lemma \ref{1}, the value \(d\) is the upper bound of \(\theta_{S(x^n)}\). Note that the sequences \(S\) and \(T\) partition \(U\). Lemma \ref{4} shows that \(P(S)\) contains a complete residue system modulo \(a\), and Lemma \ref{5} and \ref{6} together show that \(P(T)\) contains arbitrarily long arithmetic progression of integers with common difference \(a\). The properties of \(S\) and \(T\) are used in Lemma \ref{7} to prove that \(P(U)\) contains \(a\) consecutive integers.
\begin{lemma} \label{4}
The set \(P(S)\) contains a complete residue system modulo \(a\).
\begin{proof}
It suffices to prove that \(\{1,2,\ldots,a\}\subset P(\Z_a(S))\).
Let \(S_1\), \(S_2\), \ldots, \(S_{2^n}\) be the sequences defined by
\[S_j=(j^n,j^n,\ldots,j^n), \quad 1\leq j\leq 2^n\]
where \(|S_j|=2^{n^2-n}\) for all \(j\). Since for each \(0\leq \alpha\leq 2^{n^2-n}-1\), \(1\leq \beta\leq 2^n\) we have
\[f(\alpha a+\beta)\equiv \beta^n \pmod a,\]
and \(S\) is the sequence of such \(f(\alpha a+\beta)\) in increasing order, the sequences \(S_1\), \(S_2\), \ldots, \(S_{2^n}\) partition the sequence \(\Z_a(S)\). Note that \[P(S_1)=\{1,2,\ldots,2^{n^2-n}\}, \quad P(S_2)=\{2^n,2\cdot 2^n,3\cdot 2^n,\ldots,2^{n^2-n}\cdot 2^n\}.\]
Since for every integer \(1\leq m\leq 2^{n^2-n}(1+2^n)\) there exist \(0\leq\alpha\leq 2^{n^2-n}\), \(1\leq\beta\leq 2^n\) such that \[m=\alpha 2^n+\beta,\] we have \[P(S_1\cup S_2)=\{1,2,3,\ldots,2^{n^2-n}(1+2^n)\}.\]
Likewise, for every \(j\geq 3\), the inequality \[j^n<2^n(j-1)^n<2^{n^2-n}(1+2^n+\cdots+(j-1)^n)\] holds, and therefore for every \(1\leq m\leq 2^{n^2-n}(1+2^n+\cdots+j^n)\) there exists \(0\leq \alpha\leq 2^{n^2-n}\), \(1\leq \beta\leq 2^{n^2-n}(1+2^n+\cdots (j-1)^n)\) such that \(m=\alpha j^n+\beta\). Therefore
\[P(\Z_a(S))=P(S_1\cup S_2\cup\cdots\cup S_{2^n})=\{1,2,3,\ldots,2^{n^2-n}(1+2^n+3^n+\cdots+2^{n^2})\}.\]
It remains to prove that \[a=n!2^{n^2}\leq 2^{n^2-n}(1+2^n+3^n+\cdots+2^{n^2}).\]
Since \[\left(\frac{1+2^n+\cdots+2^{n^2}}{2^n}\right)^\frac{1}{n}\geq \frac{1+2+\cdots+2^n}{2^n},\] we have
\[2^{n^2-n}(1+2^n+\cdots+2^{n^2})\geq (1+2+\cdots+2^n)^n=\left(\frac{2^n(2^n+1)}{2}\right)^n.\]
Since \(2^n+1>2j\) for every positive integer \(j\leq n\), we have
\begin{align*}
    \frac{2^{n^2-n}(1+2^n+\cdots+2^{n^2})}{n!2^{n^2}}&\geq \left(\frac{2^n(2^n+1)}{2}\right)^n\cdot \frac{1}{n!2^{n^2}} \\
    &=\frac{(2^n+1)^n}{n!2^n} \\
    &=\prod_{j=1}^{n} \frac{2^n+1}{2j} \\
    &>1.
\end{align*}
Therefore, \(a=n!2^{n^2}<2^{n^2-n}(1+2^n+\cdots+2^{n^2})\) and it completes the proof.
\end{proof}
\end{lemma}
\begin{lemma} \label{5}
For every positive integer \(m\), \[a\in A\Big(\big(f(m),f(m+2),f(m+4),\ldots,f(m+\frac{2}{3}(2^{2n}-1)\big)\Big).\]
\begin{proof}
Define \(\Delta_k:\Q[x]\to \Q[x]\) by:
\begin{align*}
    \Delta_1(g(x))&=g(4x+2)-g(4x), \\
    \Delta_k(g(x))&=\Delta_1(\Delta_{k-1}(g(x))), \quad 2\leq k\leq n,
\end{align*}
so that for \(1\leq k\leq n\), \(\Delta_k(f(x))\) is a polynomial of degree \(n-k\). 
For example,
\[\Delta_2(f(x))=\Delta_1(f(4x+2)-f(4x)) \]
\[=\Big(f(16x+10)+f(16x)\Big)-\Big(f(16x+8)+f(16x+2)\Big) \]
and
\begin{align*}
    \Delta_3((f(x))&=\Delta_1(\Delta_2(f(x))) \\
    &= \Big(f(64x+42)+f(64x+32)+f(64x+8)+f(64x+2)\Big) \\
    &\quad\;\:-\Big(f(64x+40)+f(64x+34)+f(64x+10)+f(64x)\Big).
\end{align*}
It is easy to check that there are \(2^{k-1}\) positive terms and \(2^{k-1}\) negative terms in \(\Delta_k(f(x))\), and all of the terms are distinct. Therefore, for each \(1\leq k\leq n\), there exist \(2^k\) distinct integers \(\alpha_k(1)>\alpha_k(2)>\cdots>\alpha_k(2^{k-1})\), \(\beta_k(1)> \beta_k(2)>\cdots>\beta_k(2^{k-1})\) with \(\alpha_k(1)>\beta_k(1)\) such that
\[
\Delta_k(f(x))=\sum_{i=1}^{2^{k-1}}f\big(2^{2k}x+\alpha_k(i)\big)-\sum_{i=1}^{2^{k-1}}f\big(2^{2k}x+\beta_k(i)\big).
\]
Since  \(\alpha_1(1)=2\) and \(\alpha_k(1)=4\alpha_{k-1}(1)+2\) for \(k\geq 2\), we have \[\alpha_k(1)=\frac{2}{3}(2^{2k}-1).\] Also, we have \(\{\alpha_k(2^{k-1}),\beta_k(2^{k-1})\}=\{0,2\}\). Therefore
\[\Delta_k(f(x))\in A\Big(\big(f(2^{2k}x),f(2^{2k}x+2),\ldots,f(2^{2k}x+\frac{2}{3}(2^{2k}-1))\big)\Big).\] 
On the other hand, since 
\begin{align*}
    \Delta_1(f(x))&=f(4x+2)-f(4x) \\
    &=(4x+2)^n-(4x)^n \\
    &=n2^{2n-1}x^{n-1}+\textrm{terms of lower degree},
    \end{align*}
we have
\begin{align*}
     \Delta_n(f(x))&=n(n-1)(n-2)\cdots 1\cdot 2^{2n-1}2^{2n-3}2^{2n-5}\cdots 2^1 \\
     &=n!2^{n^2} \\
     &=a.
\end{align*} 
Therefore, \[a\in A\Big(\big(f(2^{2n}x),f(2^{2n}x+2),\ldots,f(2^{2n}x+\frac{2}{3}(2^{2n}-1))\big)\Big).\]
Since the \(\Delta_n(f(x))\) is a polynomial of degree \(0\), the value \(a=\Delta_n(f(x))\) is independent of \(x\). Therefore, we can replace \(2^{2n}x\) with an arbitrary positive integer \(m\) and we have 
\[
a\in A\Big(\big(f(m), f(m+2), f(m+4), \ldots, f(m+ \frac{2}{3}(2^{2n}-1))\big)\Big). \qedhere
\]
\end{proof}
\end{lemma}
\begin{lemma} \label{6}
For every positive integer \(t\), there exists a positive integer \(c\) such that all the integers \[c+ja, \quad 1\leq j\leq t\] belong to \(P(T)\) and 
\[
c<(t-1)2^{n-1}(r+\frac{2}{3}(t-1)(2^{2n}-1)+2(t-2))^n-a.
\]
\begin{proof}
Let \(\alpha=\frac{2}{3}(2^{2n}-1)\), and let \(T_1,T_2,\ldots,T_{t-1}\) be the sequences defined by
\begin{align*}
    T_1 &=\Big(f(r),f(r+2),f(r+4),\ldots,f(r+\alpha)\Big), \\
    T_2 &=\Big(f(r+\alpha+2),f(r+\alpha+4),\ldots,f(r+2\alpha+2)\Big), \\
    T_3 &=\Big(f(r+2\alpha+4),f(r+2\alpha+6),\ldots,f(r+3\alpha+4)\Big), \; \ldots \\
    T_{t-1} &=\Big(f(r+(t-2)\alpha+2(t-2)),\ldots,f(r+(t-1)\alpha+2(t-2))\Big).
\end{align*}
By Lemma \ref{5}, \(a\in A(T_j)\) for every \(1\leq j\leq t-1\), and there exists
\[A_j,B_j\in P(T_j)
\]
such that \(A_j-B_j=a\), both \(A_j\) and \(B_j\) consist of \(2^{n-1}\) terms, and all \(2^n\) terms of \(A_j\) and \(B_j\) are distinct.
Let 
\begin{align*}
C_1&=B_1+B_2+B_3+\cdots+B_{t-1}, \\
C_2&=A_1+B_2+B_3+\cdots+B_{t-1}, \\
C_3&=A_1+A_2+B_3+\cdots+B_{t-1}, \; \ldots \\
C_j&=\sum_{i=1}^{j-1}A_i+\sum_{i=j}^{t-1}B_i, \; \ldots \\
C_t&=A_1+A_2+A_3+\cdots+A_{t-1}.
\end{align*}
Then each \(C_j\) belongs to \(P(T)\), and \((C_1,C_2,\ldots,C_t)\) is an arithmetic progression of \(t\) integers with common difference \(a\). Thus, they are exactly the integers \(c+ja\), \(1\leq j \leq t\) with \(c=C_1-a=B_1+B_2+\cdots+B_{t-1}-a\). Since each \(B_j\), \(1\leq j\leq t-1\) is a sum of \(2^{n-1}\) terms in \(T\), and all of the terms are less than or equal to \[f(r+(t-1)\alpha+2(t-2))=(r+\frac{2}{3}(t-1)(2^{2n}-1)+2(t-2))^n,\] we have
\[c=C_1-a<(t-1)2^{n-1}(r+\frac{2}{3}(t-1)(2^{2n}-1)+2(t-2))^n-a.\qedhere\]
\end{proof}
\end{lemma}
Finally, we show that \(P(U)\) contains \(a\) consecutive integers \(k_1+t_1,k_2+t_2,\ldots,k_a+t_a\), where \(\{k_1,k_2,\ldots,k_a\}\) is a complete residue system of \(a\) in \(P(S)\) and \(t_1,t_2,\ldots,t_a\) are taken from the arithmetic progression in \(P(T)\).
\begin{lemma} \label{7}
Let \(b=2^{n^3}a^{n-1}\). The sequence \(U\) is \((d,a)\)-representable for a positive integer \(d\) such that
\[
d<(b-1)2^{n-1}(r+\frac{2}{3}(b-1)(2^{2n}-1)+2(b-2))^n-2a+ab.
\]
\begin{proof}
By Lemma \ref{6}, \(P(T)\) contains an arithmetic progression of \(b\) integers, \[c+ja, \quad 1\leq j\leq b\] with 
\[
c<(b-1)2^{n-1}(r+\frac{2}{3}(b-1)(2^{2n}-1)+2(b-2))^n-a.
\]
By Lemma \ref{4}, there exist positive integers \(1=k_1<k_2<\cdots<k_a\) in \(P(S)\)
such that \(\{k_1,k_2,\ldots,k_a\}\) is a complete residue system modulo \(a\). For \(1\leq j\leq a\), let
\[
n_j=\Big\lfloor\frac{k_a-k_j}{a}\Big\rfloor+1.
\]
Then for each \(1\leq j\leq a\),
\[
\frac{k_a-k_j}{a}<n_j\leq \frac{k_a-k_j}{a}+1 \iff k_a<n_ja+k_j\leq k_a+a.
\]
Also, if \(i\neq j\) then \(n_ia+k_i\not\equiv n_ja+k_j \pmod a\). Therefore 
\[
\{c+n_1a+k_1,c+n_2a+k_2,\ldots,c+n_aa+k_a\}
\]
is the set of \(a\) consecutive integers
\[
\{c+k_a+1,c+k_a+2,\ldots,c+k_a+a\}.
\]
It remains to prove that each \(c+n_ja+k_j\) is in \(P(U)\).
Let \(\Sigma(S)\) denote the sum of every element of \(S\). Since \(|S|=2^{n^2}\), and
\[
s_j\leq f((2^{n^2-n}-1)a+2^n)=(r-a+2^n)^n<r^n-(a-2^n)^n<r^n-n!
\]
for each \(s_j\in S\), we have
 \[\Sigma(S)<2^{n^2}(r^n-n!)=2^{n^2}r^n-a.\] Therefore, for each \(1\leq j\leq a\) we have
\[
1\leq n_j<\frac{k_a}{a}+1\leq\frac{1}{a}\Sigma(S)+1<\frac{1}{a}2^{n^2}r^n=2^{n^3}a^{n-1}=b
\]
and thus all of \(c+n_ja\) belong to \(P(T)\). Since all of \(k_j\) belong to \(P(S)\), all of \(c+n_ja+k_j\) belong to \(P(U)\). Therefore, \(U\) is \(
(c+k_a,a)\)-representable. Let \[d=c+k_a.\]
Since \(k_a<\Sigma(S)<2^{n^2}r^n-a=ab-a\),
\[
d=c+k_a<(b-1)2^{n-1}(r+\frac{2}{3}(b-1)(2^{2n}-1)+2(b-2))^n-2a+ab. \qedhere
\]
\end{proof}
\end{lemma}
Now we have everything we need to prove the theorem.
\section{Proof of the theorem}
\begin{proof}[\unskip\nopunct]
Recall that \(U\) and \(W\) are disjoint subsequences of \(S(f)\). By Corollary \ref{cor1}, \(W\) is a \(\Sigma(a)\)-sequence and by Lemma \ref{7}, \(U\) is \((d,a)\)-representable with 
\[
d<(b-1)2^{n-1}(r+\frac{2}{3}(b-1)(2^{2n}-1)+2(b-2))^n-2a+ab.
\]
Therefore by Lemma \ref{1}, \(S(x^n)=S(f)\) is complete and
\[
\theta_{S(x^n)}\leq d<(b-1)2^{n-1}(r+\frac{2}{3}(b-1)(2^{2n}-1)+2(b-2))^n-2a+ab. \qedhere
\]
\end{proof}
\section{Acknowledgments}
The author would like to thank Dr. Luke Oeding of Auburn University for his advice. His suggestions were valuable and helped the author to obtain a better upper bound. Also, the author would like to thank Dr. Peter Johnson of Auburn University and the anonymous referees for their helpful comments.

\bigskip
\hrule
\bigskip

\noindent 2010 {\it Mathematics Subject Classification}:  Primary
11P05; Secondary 05A17.

\noindent \emph{Keywords: } complete sequence, threshold of completeness,
sum of powers.

\bigskip
\hrule
\bigskip

\noindent (Concerned with sequence
\seqnum{A001661}.)
\end{document}